\documentclass{article}
\title{Selection principles and countable dimension} 
\author{Liljana Babinkostova and Marion Scheepers}
\date{}
\usepackage{amsfonts}
\newcommand{\naturals}{{\mathbb N}}
\newcommand{\reals}{{\mathbb R}}
\newcommand{\sone}{{\sf S}_1}
\newcommand{\gone}{{\sf G}_1}

\newtheorem{theorem}{{\bf Theorem }}
\newtheorem{corollary}[theorem]{{\bf Corollary}}
\newtheorem{lemma}[theorem]{{\bf Lemma }}

\begin{document}
\maketitle
\begin{abstract} 
We consider player TWO of the game $\gone(\mathcal{A},\mathcal{B})$ when $\mathcal{A}$ and $\mathcal{B}$ are special classes of open covers of metrizable spaces. Our results give game-theoretic characterizations of the notions of a countable dimensional and of a strongly countable dimensional metric spaces.  
\end{abstract}

The selection principle $\sone(\mathcal{A},\mathcal{B})$ states: There is for each sequence $(A_n:n\in\naturals)$ of elements of $\mathcal{A}$ a corresponding sequence $(b_n:n\in\naturals)$ such that for each $n$ we have $b_n\in A_n$, and $\{b_n:n\in\naturals\}$ is an element of $\mathcal{B}$. There are many examples of this selection principle in the literature. One of the earliest examples of it is known as \emph{the Rothberger property}, $\sone(\mathcal{O},\mathcal{O})$. Here, $\mathcal{O}$ is the collection of all open covers of a topological space. 

The following game, $\gone(\mathcal{A},\mathcal{B})$, is naturally associated with $\sone(\mathcal{A},\mathcal{B})$: Players ONE and TWO play an inning per positive integer. In the $n$-th inning ONE first chooses an element $O_n$ of $\mathcal{A}$; TWO responds by choosing an element $T_n\in O_n$. A play 
\[
  O_1,\, T_1,\, O_2,\, T_2,\, \cdots, O_n,\, T_n,\, \cdots
\]
is won by TWO if $\{T_n:n\in\naturals\}$ is in $\mathcal{B}$, else ONE wins.

\begin{center}
\begin{tabular}{c}\\
TWO has a winning strategy in $\gone(\mathcal{A},\mathcal{B})$\\
$ \Downarrow$ \\
ONE has no winning strategy in $\gone(\mathcal{A},\mathcal{B})$\\
$\Downarrow$\\
$\sone(\mathcal{A},\mathcal{B})$.
\end{tabular}
\end{center}
There are many known examples of $\mathcal{A}$ and $\mathcal{B}$ where neither of these implications reverse. 

Several classes of open covers of spaces have been defined by the following schema: For a space $X$, and a collection $\mathcal{T}$ of subsets of $X$, an open cover $\mathcal{U}$ of $X$ is said to be a $\mathcal{T}$ - cover if $X$ is not a member of $\mathcal{U}$, but there is for each $T\in\mathcal{T}$ a $U\in\mathcal{U}$ with $T\subseteq U$. The symbol $\mathcal{O}(\mathcal{T})$ denotes the collection of $\mathcal{T}$-covers of $X$. In this paper we consider only $\mathcal{A}$ which are of the form $\mathcal{O}(\mathcal{T})$ and $\mathcal{B} = \mathcal{O}$. Several examples of open covers of the form $\mathcal{O}(\mathcal{T})$ appear in the literature. To mention just a few: When $\mathcal{T}$ is the family of one-element subsets of $X$, $\mathcal{O}(\mathcal{T}) = \mathcal{O}$. When $\mathcal{T}$ is the family of finite subsets of $X$, then members of $\mathcal{O}(\mathcal{T})$ are called $\omega$-covers in \cite{G-N}. The symbol $\Omega$ denotes the family of $\omega$-covers of $X$. When $\mathcal{T}$ is the collection of compact subsets of $X$, then members of $\mathcal{O}(\mathcal{T})$ are called $k$-covers in \cite{k}. In \cite{k} the collection of $k$-covers is denoted $\mathcal{K}$.

Though some of our results hold for more general spaces, in this paper ``topological space" means separable metric space, and ``dimension" means Lebesgue covering dimension. We consider only infinite-dimensional separable metric spaces. By classical results of Hurewicz and Tumarkin these are separable metric spaces which cannot be represented as the union of finitely many zerodimensional subspaces.

\section{Properties of strategies of player TWO}

\begin{lemma}\label{twostrategy} Let $F$ be a strategy of TWO in the game $\gone(\mathcal{O}(\mathcal{T}),\mathcal{B})$. Then there is for each finite sequence $(\mathcal{U}_1,\cdots,\mathcal{U}_n)$ of elements of $\mathcal{O}(\mathcal{T})$, an element $C\in\mathcal{T}$ such that for each open set $U\supseteq C$ there is a $\mathcal{U}\in\mathcal{O}(\mathcal{T})$ such that $U = F(\mathcal{U}_1,\cdots,\mathcal{U}_n,\mathcal{U})$.
\end{lemma}
{\bf Proof:} 
For suppose on the contrary this is false. Fix a finite sequence $(\mathcal{U}_1,\cdots,\mathcal{U}_n)$ witnessing this, and choose for each set $C\subset X$ which is in $\mathcal{T}$ an open set $U_C\supseteq C$ witnessing the failure of Claim 1. Then $\mathcal{U} = \{U_C:C\subset X \mbox{ and }C\in\mathcal{T}\}$ is a member of $\mathcal{O}(\mathcal{T})$, and as $F(\mathcal{U}_1,\cdots,\mathcal{U}_n,\mathcal{U}) = U_C$ for some $C\in \mathcal{T}$, this contradicts the selection of $U_C$.
$\diamondsuit$

When $\mathcal{T}$ has additional properties, Lemma \ref{twostrategy} can be extended to reflect that. For example: The family $\mathcal{T}$ is \emph{up-directed} if there is for each $A$ and $B$ in $\mathcal{T}$, a $C$ in $\mathcal{T}$ with $A\cup B\subseteq C$. 
\begin{lemma}\label{twostrategyup} Let $\mathcal{T}$ be an up-directed family. Let $F$ be a strategy of TWO in the game $\gone(\mathcal{O}(\mathcal{T}),\mathcal{B})$. Then there is for each $D\in\mathcal{T}$ and each finite sequence $(\mathcal{U}_1,\cdots,\mathcal{U}_n)$ of elements of $\mathcal{O}(\mathcal{T})$, an element $C\in\mathcal{T}$ such that $D\subseteq C$ and for each open set $U\supseteq C$ there is a $\mathcal{U}\in\mathcal{O}(\mathcal{T})$ such that $U = F(\mathcal{U}_1,\cdots,\mathcal{U}_n,\mathcal{U})$.
\end{lemma}
{\bf Proof:} 
For suppose on the contrary this is false. Fix a finite sequence $(\mathcal{U}_1,\cdots,\mathcal{U}_n)$ and a set $D\in\mathcal{T}$ witnessing this, and choose for each set $C\subset X$ which is in $\mathcal{T}$ and with $D\subset C$ an open set $U_C\supseteq C$ witnessing the failure of Claim 1. Then, as $\mathcal{T}$ is up-directed, $\mathcal{U} = \{U_C:D\subset C\subset X \mbox{ and }C\in\mathcal{T}\}$ is a member of $\mathcal{O}(\mathcal{T})$, and as $F(\mathcal{U}_1,\cdots,\mathcal{U}_n,\mathcal{U}) = U_C$ for some $C\in \mathcal{T}$, this contradicts the selection of $U_C$.
$\diamondsuit$

We shall say that $X$ is $\mathcal{T}$-\emph{first countable} if there is for each $T\in\mathcal{T}$ a sequence $(U_n:n=1,2,\cdots)$ of open sets such that for all $n$, $T\subset U_{n+1}\subset U_n$, and for each open set $U\supset T$ there is an $n$ with $U_n\subset U$. Let $\langle\mathcal{T}\rangle$ denote the subspaces which are unions of countably many elements of $\mathcal{T}$.

\begin{theorem}\label{twoskernel} If $F$ is any strategy for TWO in $\gone(\mathcal{O}(\mathcal{T}),\mathcal{O})$ and if $X$ is $\mathcal{T}$-first countable, then there is a set $S\in\langle\mathcal{T}\rangle$ such that:
For any closed set $C\subset X\setminus S$, there is an $F$-play $O_1,\, T_1,\, \cdots,\, O_n,\, T_n\,\cdots$
such that $\bigcup_{n=1}^{\infty}T_n \subseteq X\setminus C$.
\end{theorem}

More can be proved for up-directed $\mathcal{T}$:
\begin{theorem}\label{twoskerneldirected} Let $\mathcal{T}$ be up-directed. If $F$ is any strategy for TWO in $\gone(\mathcal{O}(\mathcal{T}),\mathcal{O})$ and if $X$ is $\mathcal{T}$-first countable, then there is for each set $T\in\langle\mathcal{T}\rangle$ a set $S\in\langle\mathcal{T}\rangle$ such that:
$T\subseteq S$ and for any closed set $C\subset X\setminus S$, there is an $F$-play 
\[
  O_1,\, T_1,\, \cdots,\, O_n,\, T_n\,\cdots
\]
such that $T\subseteq \bigcup_{n=1}^{\infty}T_n \subseteq X\setminus C$.
\end{theorem}
{\bf Proof:} Let $F$ be a strategy of TWO. Let $T$ be a given element of $\langle\mathcal{T}\rangle$, and write $T=\bigcup_{n=1}^{\infty}T_n$, where each $T_n$ is an element of $\mathcal{T}$. 

Starting with $T_1$ and the empty sequence of elements of $\mathcal{O}(\mathcal{T})$, apply Lemma \ref{twostrategyup} to choose an element $S_\emptyset$ of $\mathcal{T}$ such that $T_1\subset S_\emptyset$, and for each open set $U\supseteq S_\emptyset$ there is an element $\mathcal{U}\in\mathcal{O}(\mathcal{T})$ with $U = F(\mathcal{U})$. Since $X$ is $\mathcal{T}$-first countable, choose for each $n$ an open set $U_n$ such that $U_n\supset U_{n+1}$, and for each open set $U$ with $S_\emptyset \subset U$ there is an $n$ with $U_n\subset U$. Using Lemma \ref{twostrategyup}, choose for each $n$ an element $\mathcal{U}_n$ of $\mathcal{O}(\mathcal{T})$ such that $U_n = F(\mathcal{U}_n)$.

Now consider $T_2$, and for each $n$ the one-term sequence $(\mathcal{U}_n)$ of elements of $\mathcal{O}(\mathcal{T})$. Since $\mathcal{T}$ is up-directed, choose an element $T$ of $\mathcal{T}$ with $S_\emptyset\cup T_2\subset T$. Applying Lemma \ref{twostrategyup} to $T$ and $(\mathcal{U}_n)$ choose an element $S_{(n)}\in\mathcal{T}$ such that for each open set $U\supseteq S_{(n)}$ there is a $\mathcal{U}\in\mathcal{O}(\mathcal{T})$ with $U = F(\mathcal{U}_{n},\mathcal{U})$. Since $X$ is $\mathcal{T}$-first countble, choose for each $k$ an open set $U_{(n,k)}\supseteq S_{(n)}$ such that $U_{(n,k)}\supseteq U_{(n,k+1)}\supseteq S_{(n)}$, and for each open set $U\supset S_{(n)}$ there is a $k$ with $U\supset U_{(n,k)}$. Then choose for each $n$ and $k$ an element $\mathcal{U}_{(n,k)}$ of $\mathcal{O}(\mathcal{T})$ such that $U_{(n,k}) = F(\mathcal{U}_{(n)},\mathcal{U}_{(n,k)})$.

In general, fix $k$ and suppose we have chosen for each finite sequence $(n_1,\cdots,n_k)$ of positive integers, sets $S_{(n_1,\cdots,n_k)}\in\mathcal{T}$, open sets $U_{(n_1,\cdots,n_k,n)}$ and elements $\mathcal{U}_{(n_1,\cdots,n_k,n)}$ of $\mathcal{O}(\mathcal{T})$, $n<\infty$, such that:
\begin{enumerate}
\item{$T_1\cup\cdots\cup T_k\subset S_{(n_1,\cdots,n_k)}$;}
\item{$\{U_{(n_1,\cdots,n_k,n)}:n<\infty\}$ witnesses the $\mathcal{T}$-first countability of $X$ at $S_{(n_1,\cdots,n_k)}$;}
\item{$U_{(n_1,\cdots,n_k,n)} = F(\mathcal{U}_{(n_1)},\cdots,\mathcal{U}_{(n_1,\cdots,n_k)},\mathcal{U}_{(n_1,\cdots,n_k,n)})$;}
\end{enumerate}
Now consider a fixed sequence of length $k$, say $(n_1,\cdots,n_k)$. Since $\mathcal{T}$ is up-directed choose an element $T$ of $\mathcal{T}$ such that $T_{k+1}\cup S_{(n_1,\cdots,n_k)}\subset T$. For each $n$ apply Lemma \ref{twostrategyup} to $T$ and the finite sequence $(\mathcal{U}_{(n_1)},\cdots,\mathcal{U}_{(n_1,\cdots,n_k,n)})$: Choose a set $S_{(n_1,\cdots,n_k,n)}\in\mathcal{T}$ such that $T\subseteq S_{(n_1,\cdots,n_k,n)}$ and for each open set $U\supseteq S_{(n_1,\cdots,n_k,n)}$ there is a $\mathcal{U}\in\mathcal{O}(\mathcal{T})$ such that $U = F(\mathcal{U}_{(n_1)},\cdots,\mathcal{U}_{(n_1,\cdots,n_k,n)},\mathcal{U})$. Since $X$ is $\mathcal{T}$-first countable, choose for each $j$ an open set $U_{(n_1,\cdots,n_k,n,j)}$ such that $U_{(n_1,\cdots,n_k,j+1)}\subset U_{(n_1,\cdots,n_k,n,j)}$, and for each open set $U\supset S_{(n_1,\cdots,n_k,n)}$ there is a $j$ with $U\supseteq U_{(n_1,\cdots,n_k,j)}$. Then choose for each $j$ an $\mathcal{U}_{(n_1,\cdots,n_k,n,j)}\in\mathcal{O}(\mathcal{T})$ such that $U_{(n_1,\cdots,n_k,n,j)} = F(\mathcal{U}_{(n_1)},\cdots,\mathcal{U}_{(n_1,\cdots,n_k,n)},\mathcal{U}_{(n_1,\cdots,n_k,n,j)})$.

This shows how to continue for all $k$ the recursive definition of the items $S_{(n_1,\cdots,n_k)}\in\mathcal{T}$, open sets $U_{(n_1,\cdots,n_k,n)}$ and elements $\mathcal{U}_{(n_1,\cdots,n_k,n)}$ of $\mathcal{O}(\mathcal{T})$, $n<\infty$ as above.

Finally, put $S = \cup_{\tau\in\,^{<\omega}\naturals}S_{\tau}$. It is clear that $S\in\langle\mathcal{T}\rangle$, and that $T\subset S$. Consider a closed set $C\subset X\setminus S$. Since $C\cap S_{\emptyset} = \emptyset$, choose an $n_1$ so that $U_{(n_1)}\cap C = \emptyset$. Then since $C\cap S_{(n_1)} = \emptyset$, choose an $n_2$ such that $U_{(n_1,n_2)}\cap C=\emptyset$. Since $C\cap S_{(n_1,n_2)} = \emptyset$ choose an $n_3$ so that $U_{(n_1,n_2,n_3)}\cap C = \emptyset$, and so on. In this way we find an $F$-play
\[
  \mathcal{U}_{(n_1)}, U_{(n_1)}, \mathcal{U}_{(n_1,n_2)},\, U_{(n_1,n_2)},\cdots
\] 
such that $T\subset \bigcup_{k=1}^{\infty}U_{(n_1,\cdots,n_k)} \subset X\setminus C$. $\diamondsuit$

When $\mathcal{T}$ is a collection of compact sets in a metrizable space $X$ then $X$ is $\mathcal{T}$-first countable. Call a subset $\mathcal{C}$ of $\mathcal{T}$ \emph{cofinal} if there is for each $T\in\mathcal{T}$ a $C\in\mathcal{C}$ with $T\subseteq C$. As an examination of the proof of Theorem \ref{twoskerneldirected} reveals, we do not need full $\mathcal{T}$-first countability of $X$, but only that $X$ is $\mathcal{C}$-first countable for some cofinal set $\mathcal{C}\subseteq\mathcal{T}$. Thus, we in fact have:

\begin{theorem}\label{twoskerneldirectedcofinal} Let $\mathcal{T}$ be up-directed. If $F$ is any strategy for TWO in $\gone(\mathcal{O}(\mathcal{T}),\mathcal{O})$ and if $X$ is $\mathcal{C}$-first countable where $\mathcal{C}\subset\mathcal{T}$ is cofinal in $\mathcal{T}$, then there is for each set $T\in\langle\mathcal{T}\rangle$ a set $S\in\langle\mathcal{C}\rangle$ such that:
$T\subseteq S$ and for any closed set $C\subset X\setminus S$, there is an $F$-play 
\[
  O_1,\, T_1,\, \cdots,\, O_n,\, T_n\,\cdots
\]
such that $T\subseteq \bigcup_{n=1}^{\infty}T_n \subseteq X\setminus C$.
\end{theorem}

\section{When player TWO has a winning strategy}

Recall that a subset of a topological space is a ${\sf G}_{\delta}$-set if it is an intersection of countably many open sets.

\begin{theorem}\label{generalplayertwo} If the family $\mathcal{T}$ has a cofinal subset consisting of ${\sf G}_{\delta}$ subsets of $X$, then TWO has a winning strategy in $\gone(\mathcal{O}(\mathcal{T}),\mathcal{O})$ if, and only if, the space is a union of countably many members of $\mathcal{T}$.
\end{theorem}
{\bf Proof:} $2\Rightarrow 1$ is easy to prove. We prove $1\Rightarrow 2$. Let $F$ be a winning strategy for TWO. Let $\mathcal{C}\subseteq \mathcal{T}$ be a cofinal set consisting of ${\sf G}_{\delta}$-sets.\\
By Lemma \ref{twostrategy} choose $C_{\emptyset}\in\mathcal{T}$ associated to the empty sequence. Since $\mathcal{C}$ is cofinal in $\mathcal{T}$, choose for $C_{\emptyset}$ a ${\sf G}_{\delta}$ set $G_{\emptyset}$ in $\mathcal{C}$ with $C_{\emptyset}\subseteq G_{\emptyset}$. Choose open sets $(U_n:n\in\naturals)$ such that for each $n$ we have $G_{\emptyset}\subset U_{n+1}\subset U_n$, and $G_{\emptyset} = \cap_{n\in\naturals}U_n$.   

For each $n$ choose by Lemma \ref{twostrategy} a cover $\mathcal{U}_n\in\mathcal{O}({\mathcal{T}})$ with $U_n = F(\mathcal{U}_n)$. Choose for each $n$ a  $C_n\in\mathcal{T}$ associated to $(\mathcal{U}_n)$ by Lemma \ref{twostrategy}. For each $n$ also choose a ${\sf G}_{\delta}$-set $G_n\in\mathcal{C}$ with $C_n\subseteq G_n$. For each $n_1$ choose a sequence $(U_{n_1n}:n\in\naturals)$ of open sets such that $G_{n_1} = \cap_{n\in\naturals}U_{n_1n}$ and for each $n$, $U_{n_1n+1}\subset U_{n_1n}$. For each $n_1n_2$ choose by Lemma \ref{twostrategy} a cover $\mathcal{U}_{n_1n_2}\in\mathcal{O}(\mathcal{T})$ such that $U_{n_1n_2} = F(\mathcal{U}_{n_1},\mathcal{U}_{n_1n_2})$. Choose by Lemma \ref{twostrategy} a $C_{n_1n_2}\in\mathcal{T}$ associated to $(\mathcal{U}_{n_1},\mathcal{U}_{n_1n_2})$, and then choose a ${\sf G}_{\delta}$-set $G_{n_1n_2}\in\mathcal{C}$ with $C_{n_1n_2}\subset G_{n_1n_2}$, and so on. 

Thus we get for each finite sequence $(n_1n_2\cdots n_k)$ of positive integers 
\begin{enumerate}
\item{a set $C_{n_1 \cdots n_k}\in\mathcal{T}$,}
\item{a ${\sf G}_{\delta}$-set $G_{n_1 \cdots n_k}\in\mathcal{T}$ with $C_{n_1\cdots n_k}\subseteq G_{n_1\cdots n_k}$,}
\item{a sequence $(U_{n_1\cdots n_kn}:n\in\naturals)$ of open sets with  $G_{n_1 \cdots n_k} = \cap_{n\in\naturals}U_{n_1\cdots n_kn}$ and for each $n$ $U_{n_1 \cdots n_kn+1}\subseteq U_{n_1 \cdots n_kn}$, and} 
\item{a $\mathcal{U}_{n_1\cdots n_k}\in\mathcal{O}_{(\mathcal{T})}$ such that for all $n$
\[
  U_{n_1\cdots n_kn} = F(\mathcal{U}_{n_1},\cdots,\mathcal{U}_{n_1\cdots n_kn}).
\] 
}
\end{enumerate}

Now $X$ is the union of the countably many sets $G_{\tau}\in\mathcal{T}$ where $\tau$ ranges over $^{<\omega}\, \naturals$. For if not, choose $x\in X$ which is not in any of these sets. Since $x$ is not in $G_{\emptyset}$, choose $U_{n_1}$ with $x\not\in U_{n_1}$. Now $x$ is not in $G_{n_1}$, so choose $U_{n_1n_2}$ with $x\not\in U_{n_1n_2}$, and so on. In this way we obtain the $F$-play
\[
  \mathcal{U}_{n_1},\, U_{n_1},\, \mathcal{U}_{n_1n_2},\, U_{n_1n_2},\, \cdots 
\]
lost by TWO, contradicting that $F$ is a winning strategy for TWO.
$\diamondsuit$

Examples of up-directed families $\mathcal{T}$ include:
\begin{itemize}
\item{$[X]^{<\aleph_0}$, the collection of finite subsets of $X$;}
\item{$\mathcal{K}$, the collection of compact subsets of $X$;}
\item{${\sf KFD}$, the collection of compact, finite dimensional subsets of $X$.}
\item{${\sf CFD}$, the collection of closed, finite dimensional subsets of $X$.}
\item{${\sf FD}$, the collection of finite dimensional subsets of $X$.}
\end{itemize}

A subset of a topological space is said to be \emph{countable dimensional} if it is a union of countably many zero-dimensional subsets of the space. A subset of a space is \emph{strongly countable dimensional} if it is a union of countably many closed, finite dimensional subsets. Let $X$ be a space which is not finite dimensional. Let $\mathcal{O}_{\sf cfd}$ denote $\mathcal{O}({\sf CFD})$, the collection of {\sf CFD}-covers of $X$. And let $\mathcal{O}_{{\sf fd}}$ denote $\mathcal{O}({\sf FD})$, the collection of ${\sf FD}$-covers of $X$.

\begin{corollary}\label{playerIIscd} For a metrizable space $X$ the following are equivalent:
\begin{enumerate}
\item{$X$ is strongly countable dimensional.}
\item{TWO has a winning strategy in $\gone(\mathcal{O}_{\sf cfd},\mathcal{O})$.}
\end{enumerate}
\end{corollary}
{\bf Proof:} $1\Rightarrow 2$ is easy to prove. To see $2\Rightarrow 1$, observe that in a metric space each closed set is a ${\sf G}_{\delta}$-set. Thus, $\mathcal{T}={\sf CFD}$ meets the requirements of Theorem \ref{generalplayertwo}.
$\diamondsuit$

For the next application we use the following classical theorem of Tumarkin:

\begin{theorem}[Tumarkin]\label{tumarkin} In a separable metric space each $n$-dimensional set is contained in an $n$-dimensional ${\sf G}_{\delta}$-set. 
\end{theorem}

\begin{corollary}\label{playerIIcd} For a separable metrizable space $X$ the following are equivalent:
\begin{enumerate}
\item{$X$ is countable dimensional.}
\item{TWO has a winning strategy in $\gone(\mathcal{O}_{\sf fd},\mathcal{O})$.}
\end{enumerate}
\end{corollary}
{\bf Proof:} $1\Rightarrow 2$ is easy to prove. We now prove $2\Rightarrow 1$. By Tumarkin's Theorem, $\mathcal{T} = {\sf FD}$ has a cofinal subset consisting of ${\sf G}_{\delta}$-sets. Thus the requirements of Theorem \ref{generalplayertwo} are met.
$\diamondsuit$

Recall that a topological space is \emph{perfect} if every closed set is a ${\sf G}_{\delta}$-set.

\begin{corollary}\label{kcovers} In a perfect space the following are equivalent:
\begin{enumerate}
\item{TWO has a winning strategy in $\gone(\mathcal{K},\mathcal{O})$.}
\item{The space is $\sigma$-compact.}
\end{enumerate}
\end{corollary}
{\bf Proof:} In a perfect space the collection of closed sets are ${\sf G}_{\delta}$-sets. Apply Theorem \ref{generalplayertwo}. $\diamondsuit$

And when $\mathcal{T}$ is up-directed, Theorem \ref{generalplayertwo} can be further extended to:

\begin{theorem}\label{updirectedII} If $\mathcal{T}$ is up-directed and has a cofinal subset consisting of ${\sf G}_{\delta}$-subsets of $X$, the following  are equivalent:
\begin{enumerate}
\item{TWO has a winning strategy in $\gone(\mathcal{O}(\mathcal{T}),\Gamma)$.}
\item{TWO has a winning strategy in $\gone(\mathcal{O}(\mathcal{T}),\Omega)$.}
\item{TWO has a winning strategy in $\gone(\mathcal{O}(\mathcal{T}),\mathcal{O})$.}
\end{enumerate}
\end{theorem}
{\bf Proof:} We must show that $3\Rightarrow 1$. Since $X$ is a union of countably many sets in $\mathcal{T}$, and since $\mathcal{T}$ is up-directed, we may represent $X$ as $\bigcup_{n=1}^{\infty}X_n$ where for each $n$ we have $X_n\subset X_{n+1}$ and $X_n\in\mathcal{T}$. Now, when ONE presents TWO with $O_n\in\mathcal{O}(\mathcal{T})$ in inning $n$, then TWO chooses $T_n\in O_n$ with $X_n\subset T_n$. The sequence of $T_n$'s chosen by TWO in this way results in a $\gamma$-cover of $X$. $\diamondsuit$

\section{Longer games and player TWO}

Fix an ordinal $\alpha$. Then the game $\gone^{\alpha}(\mathcal{A},\mathcal{B})$ has $\alpha$ innings and is played as follows. In inning $\beta$ ONE first chooses an $O_{\beta}\in\mathcal{A}$, and then TWO responds with a $T_{\beta}\in O_{\beta}$. A play 
\[
 O_0, T_0,\cdots,O_{\beta},  T_{\beta},\cdots,\, \beta<\alpha
\]
is won by TWO if $\{T_{\beta}:\beta<\alpha\}$ is in $\mathcal{B}$; else, ONE wins.

In this notation the game $\gone(\mathcal{A},\mathcal{B})$ is $\gone^{\omega}(\mathcal{A},\mathcal{B})$. For a space $X$ and a family $\mathcal{T}$ of subsets of $X$ with $\cup\mathcal{T} = X$, define: 
\[
  {\sf cov}_X(\mathcal{T})=\min\{|\mathcal{S}|:\mathcal{S}\subseteq \mathcal{T} \mbox{ and } X = \cup\mathcal{S}\}.
\]
When $X = \cup\mathcal{T}$, there is an ordinal $\alpha\le {\sf{cov}}_X(\mathcal{T})$ such that TWO has a winning strategy in $\gone^{\alpha}(\mathcal{O}(\mathcal{T}),\mathcal{O})$. In general, there is an ordinal $\alpha\le |X|$ such that TWO has a winning strategy in $\gone^{\alpha}(\mathcal{O}(\mathcal{T}),\mathcal{O})$.
\[
  {\sf tp}_{\sone(\mathcal{O}(\mathcal{T}),\mathcal{O})}(X) = \min\{\alpha: \mbox{ TWO has a winning strategy in }\gone^{\alpha}(\mathcal{O}(\mathcal{T}),\mathcal{O})\}.
\] 

\subsection{General properties}

The proofs of the general facts in the following lemma are left to the reader.
\begin{lemma}\label{properties}
\begin{enumerate}
\item{If $Y$ is a closed subset of $X$ then $  {\sf tp}_{\sone(\mathcal{O}(\mathcal{T}),\mathcal{O})}(Y)\le   {\sf tp}_{\sone(\mathcal{O}(\mathcal{T}),\mathcal{O})}(X)$.}
\item{If $\alpha$ is a limit ordinal and if $  {\sf tp}_{\sone(\mathcal{O}(\mathcal{T}),\mathcal{O})}(X_n)\le \alpha$ for each $n$, then $  {\sf tp}_{\sone(\mathcal{O}(\mathcal{T}),\mathcal{O})}(\bigcup_{n<\infty}X_n)\le \alpha$.}
\end{enumerate}
\end{lemma}

We shall now give examples of ordinals $\alpha$ for which TWO has winning strategies in games of length $\alpha$. First we have the following general lemma. 

\begin{lemma}\label{ordinallengthtwo} Let $X$ be $\mathcal{T}$-first countable. Assume that:
\begin{enumerate}
\item{$\mathcal{T}$ is up-directed;}
\item{$X\not\in\langle\mathcal{T}\rangle$;}
\item{$\alpha$ is the least ordinal such that there is an element $B$ of $\langle\mathcal{T}\rangle$ such that for any closed set $C\subset X\setminus B$ with $C\not\in\mathcal{T}$,  ${\sf tp}_{\sone(\mathcal{O}(\mathcal{T}),\mathcal{O})}(C)\le\alpha$.}
\end{enumerate}
 Then ${\sf tp}_{\sone(\mathcal{O}(\mathcal{T}),\mathcal{O})}(X) = \omega+\alpha$.
\end{lemma}
{\bf Proof:} We must show that TWO has a winning strategy for $\gone^{\omega+\alpha}(\mathcal{O}(\mathcal{T}),\mathcal{O})$, and that there is no $\beta<\omega+\alpha$ for which TWO has a winning strategy in $\gone^{\beta}(\mathcal{O}(\mathcal{T}),\mathcal{O})$. 

To see that TWO has a winning strategy in $\gone^{\omega+\alpha}(\mathcal{O}(\mathcal{T}),\mathcal{O})$, fix a $B$ as in the hypothesis, and for each closed set $F$ disjoint from $B$, fix a winning strategy $\tau_F$ for TWO in the game $\gone^{\alpha}(\mathcal{O}(\mathcal{T}),\mathcal{O})$ played on $F$. Now define a strategy $\sigma$ for TWO  in $\gone^{\omega+\alpha}(\mathcal{O}(\mathcal{T}),\mathcal{O})$ on $X$ as follows: During the first $\omega$ innings, TWO covers $B$. Let $T_1,\, T_2,\, \cdots$ be TWO's moves during these $\omega$ innings, and put $C = X\setminus \bigcup_{n=1}^{\infty}T_n$. Then $C$ is a closed subset of $X$, disjoint from $B$. Now TWO follows the strategy $\tau_C$ in the remaining $\alpha$ innings, to also cover $C$. 

To see that there is no $\beta<\omega+\alpha$ for which TWO has a winning strategy in $\gone^{\beta}(\mathcal{O}(\mathcal{T}),\mathcal{O})$, argue as follows: Suppose on the contrary that $\beta<\omega+\alpha$ is such that TWO has a winning strategy $\sigma$ for $\gone^{\beta}(\mathcal{O}(\mathcal{T}),\mathcal{O})$ on $X$. We will show that there is a set $S\in\langle\mathcal{T}\rangle$ and an ordinal $\gamma<\alpha$ such that for each closed set $C$ disjoint from $S$, TWO has a winning strategy in $\gone^{\gamma}(\mathcal{O}(\mathcal{T}),\mathcal{O})$ on $C$. This gives a contradiction to the minimality of $\alpha$ in hypothesis 3.

We consider cases: First, it is clear that $\alpha\le \beta$, for otherwise TWO may merely follow the winning strategy on $X$ and relativize to any closed set $C$ to win on $C$ in $\beta<\alpha$ innings, a contradiction. Thus, $\omega+\alpha > \alpha$. Then we have $\alpha<\omega^2$, say $\alpha = \omega\cdot n + k$. Since then $\omega+\alpha = \omega\cdot(n+1) + k$, we have that $\beta$ with $\alpha \le \beta <\omega+\alpha$ has the form $\beta = \omega\cdot n +\ell$ with $\ell \ge k$. The other possibility, $\beta = \omega\cdot(n+1)+j$ for some $j<k$, does not occur because it would give
$\alpha+\omega > \beta = \omega\cdot n + (\omega+j) = (\omega\cdot n +k)+(\omega+j) =\alpha+\omega+j$.

Let $F$ be a winning strategy for TWO in $\gone^{\beta}(\mathcal{O}(\mathcal{T}),\mathcal{O})$. By the second hypothesis and Theorem \ref{generalplayertwo} we have $\beta>\omega$. By Theorem \ref{twoskerneldirected} fix an element $S\in\langle T\rangle$ such that $B\subset S$, and for any closed set $C\subset X\setminus S$, there is an $F$-play $(O_1,\, T_1,\, \cdots,O_n,\, T_n,\,\cdots)$ with $S\subset (\bigcup_{n=1}^{\infty}T_n)$, and $C\cap (\bigcup_{n=1}^{\infty}T_n) = \emptyset$. Choose a closed set $C \subset X\setminus S$ with $C\not\in\mathcal{T}$. This is possible by the second hypothesis. Choose an $F$-play $(O_1,\, T_1,\, \cdots,O_n,\, T_n,\,\cdots)$ with $S\subset (\bigcup_{n=1}^{\infty}T_n)$, and $C\cap (\bigcup_{n=1}^{\infty}T_n) = \emptyset$. This $F$-play contains the first $\omega$ moves of a play according to the winning strategy $F$ for TWO in $\gone^{\beta}(\mathcal{O}(\mathcal{T}),\mathcal{O})$, and using it as strategy to play this game on $C$, we see that it requires (an additional) $\gamma = \omega\cdot(n-1)+\ell < \alpha$ innings for TWO to win on $C$. Here, $\ell$ is fixed and the same for all such $C$. Thus: ${\sf tp}_{\sone(\mathcal{O}(\mathcal{T}),\mathcal{O})}(C) \le \gamma < \alpha$. This is in contradiction to the minimality of $\alpha$. 
$\diamondsuit$

\subsection{Examples}

For each $n$ put $\reals_n = \{x\in \reals^{\naturals}:(\forall m>n)(x(m)=0)\}$. Then $\reals_n$ is homeomorphic to $\reals^n$ and thus is $\sigma$-compact, and $n$-dimensional. Thus $\reals_{\infty} = \bigcup_{n=1}^{\infty}\reals_n$ is a $\sigma$-compact strongly countable dimensional subset of $\reals^{\naturals}$. 

We shall now use the Continuum Hypothesis to construct for various infinite countable ordinals $\alpha$ subsets of $\reals^{\naturals}$ in which TWO has a winning strategy in $\gone^{\alpha}(\mathcal{O}(\mathcal{T}),\mathcal{O})$. The following is one of our main tools for these constructions:

\begin{lemma}\label{gdeltas} If $G$ is any ${\sf G}_{\delta}$-subset of $\reals^{\naturals}$ with $\reals_{\infty}\subset G$, then $G\setminus\reals_{\infty}$ contains a compact nowhere dense subset $C$ which is homeomorphic to $[0,1]^{\naturals}$.
\end{lemma}

We call $[0,1]^{\naturals}$ the Hilbert cube. From now on assume the Continuum Hypothesis. Let $(F_{\alpha}:\alpha<\omega_1)$ enumerate all the finite dimensional ${\sf G}_{\delta}$-subsets of $\reals^{\naturals}$, and let $(C_{\alpha}:\alpha<\omega_1)$ enumerate the ${\sf G}_{\delta}$-subsets which contain $\reals_{\infty}$. Recursively choose compact sets $D_{\alpha} \subset \reals^{\naturals}$, each homeomorphic to the Hilbert cube and nowhere dense, such that $D_0\subset C_{0}\setminus(\reals_{\infty}\cup F_0)$, and for all $\alpha>0$,
\[
  D_{\alpha} \subset (\cap_{\beta\le\alpha}C_{\beta})\setminus(\reals_{\infty}\cup(\bigcup\{D_{\beta}:\beta<\alpha\})\cup(\bigcup_{\beta\le\alpha}F_{\beta})).
\]
{\flushleft{\bf Version 1:}}
For each $\alpha$, choose a point $x_{\alpha}\in D_{\alpha}$ and put 
\[
  B:=\reals_{\infty}\cup\{x_{\alpha}:\alpha<\omega_1\}.
\]
{\flushleft{\bf Version 2:}}
For each $\alpha$, choose a strongly countable dimensional set $S_{\alpha}\subset D_{\alpha}$ and put 
\[
  B:=\reals_{\infty}\cup(\bigcup\{S_{\alpha}:\alpha<\omega_1\}).
\]
{\flushleft{\bf Version 3:}}
For each $\alpha$, choose a countable dimensional set $S_{\alpha}\subset D_{\alpha}$ and put 
\[
  B:=\reals_{\infty}\cup(\bigcup\{S_{\alpha}:\alpha<\omega_1\}).
\]

In all three versions, $B$ is not countable dimensional: Otherwise it would be, by Tumarkin's Theorem, for some $\alpha<\omega_1$ a subset of $\bigcup_{\beta<\alpha}F_{\beta}$. Thus TWO has no winning strategy in the games $\gone(\mathcal{O}_{{\sf cfd}},\mathcal{O})$ and $\gone(\mathcal{O}_{{\sf fd}},\mathcal{O})$. Also, in all three versions the elements of the family $\mathcal{C}$ of finite unions of the sets $S_{\alpha}$ are ${\sf G}_{\delta}$-sets in $X$, and in fact $X$ is $\mathcal{C}$-first-countable. This is because the $D_{\alpha}$'s are compact and disjoint, and $\reals^{\naturals}$ is $\mathcal{D}$-first countable, where $\mathcal{D}$ is the family of finite unions of the $D_{\alpha}$'s, and this relativizes to $X$. 

For Version 1 TWO has a winning strategy in $\gone^{\omega+1}(\mathcal{O}_{{\sf cfd}},\mathcal{O})$ and in $\gone^{\omega+1}(\mathcal{O}_{{\sf fd}},\mathcal{O})$, and in $\gone^{\omega+\omega}(\mathcal{K},\mathcal{O})$. For Version 2 TWO has a winning strategy in $\gone^{\omega+\omega}(\mathcal{O}_{{\sf cfd}},\mathcal{O})$, and for Version 3 TWO has a winning strategy in $\gone^{\omega+\omega}(\mathcal{O}_{{\sf fd}},\mathcal{O})$. 

To see this, note that in the first $\omega$ innings, TWO covers $\reals_{\infty}$. Let $\{U_n:n\in\naturals\}$ be TWO's responses in these innings. Then $G = \bigcup_{n=1}^{\infty}U_n$ is an open set containing $\reals_{\infty}$, and so there is an $\alpha<\omega_1$ such that:
{\flushleft{\bf Version 1:}}  $B\setminus G \subseteq \{x_{\beta}:\beta<\alpha\}$ is a closed, countable subset of $X$ and thus closed, zero-dimensional. In inning $\omega+1$ TWO chooses from ONE's cover an element containing the set $B\setminus G$.\\ 
{\flushleft{\bf Version 2:}}
$B\setminus G \subseteq \bigcup_{\beta<\alpha}S_{\beta}$. But $\bigcup_{\beta<\alpha}S_{\alpha}$ is strongly countable dimensional, and so TWO can cover this part of $B$ in the remaining $\omega$ innings. By Lemma \ref{ordinallengthtwo} TWO does not have a winning strategy in fewer then $\omega+\omega$ innings. 
{\flushleft{\bf Version 3:}}
$B\setminus G \subseteq \bigcup_{\beta<\alpha}S_{\beta}$. But $\bigcup_{\beta<\alpha}S_{\alpha}$ is strongly countable dimensional, and so TWO can cover this part of $B$ in the remaining $\omega$ innings. By Lemma \ref{ordinallengthtwo} TWO does not have a winning strategy in fewer then $\omega+\omega$ innings. 

With these examples established, we can now upgrade the construction as follows: Let $\alpha$ be a countable ordinal for which we have constructed an example of a subspace $S$ of $\reals^{\naturals}$ for which ${\sf tp}_{\sone(\mathcal{O}(\mathcal{T}),\mathcal{O})}(S) = \alpha$. Then choose inside each $D_{\beta}$ a set $C_{\beta}$ for which ${\sf tp}_{\sone(\mathcal{O}(\mathcal{T}),\mathcal{O})}(C_{\beta}) = \alpha$. Then the resulting subset $B$ constructed above has, by Lemma \ref{ordinallengthtwo},  ${\sf tp}_{\sone(\mathcal{O}(\mathcal{T}),\mathcal{O})}(B) = \omega+\alpha$. In this way we obtain examples for each of the lengths $\omega\cdot n$ and $\omega\cdot n+1$, for all finite $n$.

By taking topological sums and using part 2 of Lemma \ref{properties} we get examples for $\omega^2$.

\section{Conclusion}

One obvious question is whether there is, under the Continuum Hypothesis, for each limit ordinal $\alpha$ subsets $X_{\alpha}$ and $Y_{\alpha}$ of $\reals^{\naturals}$ such that ${\sf tp}_{\sone(\mathcal{O}_{{\sf cfd}},\mathcal{O})}(X_{\alpha}) = \alpha$, and ${\sf tp}_{\sone(\mathcal{O}_{{\sf cfd}},\mathcal{O})}(Y_{\alpha}) = \alpha+1$. And the same question can be asked for ${\sf tp}_{\sone(\mathcal{O}_{{\sf fd}},\mathcal{O})}$.

In \cite{babinkostova1} countable dimensionality of metrizable spaces were characterized in terms of the selective screenability game. A natural question is how $\sone(\mathcal{O}_{{\sf fd}},\mathcal{O})$ and $\sone(\mathcal{O}_{{\sf cfd}},\mathcal{O})$ are related to selective screenability. It is clear that $\sone(\mathcal{O}_{{\sf fd}},\mathcal{O})\Rightarrow \sone(\mathcal{O}_{{\sf cfd}},\mathcal{O})$. The relationship among these two classes and selective screenability is further investigated in \cite{babinkostova2} where it is shown, for example, that $\sone(\mathcal{O}_{{\sf cfd}},\mathcal{O})$ implies selective screenability, but the converse does not hold. Thus, these two classes are new classes of weakly infinite dimensional spaces.

\bigskip

\begin{center}
  Addresses
\end{center}

\medskip

\begin{flushleft}
Liljana Babinkostova                      \\
Boise State University                    \\
Department of Mathematics                 \\
Boise, ID 83725 USA                       \\
e-mail: liljanab@math.boisestate.edu      \\
\end{flushleft}

\begin{flushleft}
Marion Scheepers                    \\
Department of Mathematics           \\
Boise State University              \\
Boise, Idaho 83725 USA              \\
e-mail: marion@math.boisestate.edu   \\
\end{flushleft}

\end{document}